\documentclass{article}

\usepackage[english]{babel}
\usepackage[T1]{fontenc}
\usepackage{lmodern}
\usepackage[margin=2.69506cm]{geometry}

\usepackage{amsmath,amssymb,amsthm}
\usepackage{graphicx,subcaption,float}
\usepackage[export]{adjustbox} 
\usepackage{comment}
\usepackage{pgfplots} 
\pgfplotsset{compat=1.18}
\usepackage{microtype}
\usepackage{hyperref}
\usepackage[backend=biber,style=alphabetic,sorting=nyt,maxnames=99,doi=false,url=false,isbn=false]{biblatex}
\newcommand{\perm}{\mathrm{perm}}
\newcommand{\Hpn}{H^{\perm}_n}
\newcommand{\Hp}{H^{\perm}}

\newtheorem{theorem}{Theorem}[section]
\newtheorem{corollary}[theorem]{Corollary}
\newtheorem{lemma}[theorem]{Lemma}
\newtheorem{proposition}[theorem]{Proposition}
\theoremstyle{definition}
\newtheorem{definition}[theorem]{Definition}
\newtheorem{remark}[theorem]{Remark}

\title{The generic behavior of permutation entropy for continuous maps}
\author{Balázs Maga\thanks{HUN-REN Alfréd Rényi Institute of Mathematics, Reáltanoda utca 13-15, 1053 Budapest. E-mail: {\tt magab@renyi.hu}. Supported by the
NKFIH Grant no. 152535, 
funded by the Ministry of Innovation and Technology of Hungary from the National Research, Development and Innovation Fund.} \and Sára Szepessy\thanks{ELTE Eötvös Loránd University, Pázmány Péter Sétány 1/C, 1117 Budapest. E-mail: {\tt sara.szepessy@gmail.com}}}
\date{}

\begin{document}

\maketitle

\begin{abstract}
Permutation entropy, introduced by Bandt and Pompe in 2002, measures the complexity of a dynamical system through the ordinal patterns of its trajectories. Bandt, Keller and Pompe established that for measure-preserving, piecewise monotone interval maps with finitely many monotonicity intervals, the normalized order-$n$ permutation entropy $H_n^{\mathrm{perm}}(f)/n$ converges to the Kolmogorov--Sinai entropy. Although this result has subsequently been extended to broader classes of maps, whether normalized permutation entropy always converges has remained unknown.

We show that convergence can fail, and indeed that nonconvergence is generic. More precisely, for a generic continuous Lebesgue measure-preserving map $f:[0,1]\to[0,1]$, we prove that
\[
\liminf_{n\to\infty}\frac{H_n^{\mathrm{perm}}(f)}{n}=0,
\qquad
\limsup_{n\to\infty}\frac{H_n^{\mathrm{perm}}(f)}{n\log n}=1.
\]
Thus, despite the generic vanishing of Kolmogorov--Sinai entropy, permutation entropy exhibits oscillations between vanishing linear growth and the maximal possible asymptotic growth.

We also construct an explicit Lebesgue measure-preserving interval map with countably many monotonicity intervals for which the normalized permutation entropy has lower limit zero and upper limit infinity. This demonstrates that nonconvergence can occur even within the countably piecewise-monotone setting, beyond the scope of the convergence theorem of Gutjahr and Keller, relying on conditions concerning the lengths of the pieces.


\end{abstract}
\section{Introduction}

Introduced by Bandt and Pompe in 2002 \cite{Bandt2002}, permutation entropy is a nonparametric, ordinal-pattern-based measure of time-series complexity. For a given map $f:[0,1]\to [0,1]$, the ordinal pattern of length $n$ at a point $x$ is defined by the relative order of the sequence of iterates $x, f(x),\ldots, f^{n-1}(x)$, given that all these iterates are distinct. If periodic points form a Lebesgue null set, taking $x\sim \mathrm{Unif}(0, 1)$, this gives rise to a permutation-valued random variable. The order-$n$ permutation entropy, denoted as $\Hpn(f)$, is then defined as the Shannon entropy of this permutation-valued random variable. The \emph{permutation entropy} of $f$ is $\Hp(f)=\lim_n \Hpn(f)/n$, whenever the limit exists.

For piecewise monotone interval maps with finitely many monotonicity intervals preserving the Lebesgue measure, Bandt, Keller and Pompe established that the permutation entropy exists and it coincides with the Kolmogorov-Sinai entropy of $f$ \cite{BandtKellerPompe}. This foundational result was later generalized to broader classes of dynamical systems (see, for example, \cite[Theorem 2]{GK}, allowing countably many monotone pieces once the discrete random variable with probability masses given by the piece lengths has finite Shannon entropy). 
This regular convergence might suggest that the permutation entropy rate typically stabilizes to a single finite limit. However, the assumption of piecewise monotonicity is quite restrictive, and whether this limit necessarily exists in general is a surprisingly unexplored question. Accordingly, the literature commonly works with lower and upper permutation entropy, defined using the liminf and limsup of the normalized entropy sequence, respectively. As highlighted in a recent review by Gutjahr and Keller \cite{GK2020}, determining the conditions under which the KS entropy coincides with the upper or lower permutation entropy remains an open problem. In fact, they explicitly note that ``there is not even known an example where the entropies differ'' \cite[Section~4]{GK2020}. 

This apparent gap motivates our study. Rather than searching for a single counterexample, we investigate the generic behavior (in the sense of Baire category) of permutation entropy in the space $\mathcal{C}(\lambda)$ of continuous, Lebesgue measure-preserving functions on $[0,1]$. ($\lambda$ denotes the Lebesgue measure throughout.) Our main result shows that nonconvergence is generic in this space, and indeed occurs with the strongest possible asymptotic oscillations.

\begin{theorem} \label{mainthm}
    For a generic $f\in \mathcal{C}(\lambda)$, we have 
    $$\liminf_{n\to \infty} \frac{\Hpn (f)}{n}=0, \qquad \limsup_{n\to \infty} \frac{\Hpn (f)}{n \log n}=1.$$
\end{theorem}

\begin{remark}
    Owing to the fact that lower permutation entropy bounds the Kolmogorov--Sinai entropy from above  \cite[Theorem~1]{GK2020}, the first equality can be considered an improvement of the recent result by Bobok, \v{C}inc, Oprocha, and Troubetzkoy \cite{BCOT26}, demonstrating that the Kolmogorov--Sinai entropy is generically zero. 
\end{remark}

\begin{remark}
 The second equality demonstrates that the upper permutation entropy achieves the maximal theoretically possible growth. Since the number of possible ordinal patterns of length $n$ is $n!$, the permutation entropy is trivially bounded by $\log n! \sim n \log n$. Theorem \ref{mainthm} shows that for a generic system, the entropy sequence approximates this absolute theoretical upper bound along suitable subsequences.
\end{remark}

\begin{remark}
    A related phenomenon was recently established for sampling entropy, a different measure of permutation complexity obtained by sampling independent points from the graph of a measure-preserving map. In earlier work, the first author proved analogous generic oscillations for sampling entropy \cite[Theorem~1.4]{MB}. The present theorem shows that this phenomenon persists for temporal ordinal patterns, despite the fundamentally different sampling mechanism.
\end{remark}

Before proving Theorem \ref{mainthm}, in Section \ref{sec:examples} we explicitly construct two piecewise monotone functions (with countably infinitely many pieces) for which $\Hpn$ behaves in a peculiar way, hence demonstrating the limitations of the aforementioned convergence theorem of Gutjahr and Keller \cite[Theorem~2]{GK}. The first construction demonstrates that the lower permutation entropy can be infinite even when the Kolmogorov--Sinai entropy vanishes:

\begin{proposition} \label{prop:1}
    There exists a Lebesgue measure-preserving map $f:[0, 1]\to [0, 1]$ which is piecewise monotone with countably many monotone pieces, has no periodic points, and for which
    $$h_{KS}(f)=0, \quad \lim_{n\to\infty}\frac{\Hpn (f)}{n}=\infty.$$
\end{proposition}

This second construction gives a piecewise monotone example where the permutation entropy fails to exist:

\begin{proposition} \label{prop:2}
    There exists a Lebesgue measure-preserving map $f:[0, 1]\to [0, 1]$ which is piecewise monotone with countably many monotone pieces, has no periodic points, and for which
    $$h_{KS}(f)=\liminf_{n\to\infty}\frac{\Hpn (f)}{n}=0, \qquad \limsup_{n\to\infty}\frac{\Hpn (f)}{n}=\infty.$$
\end{proposition}

We remark that by tuning the parameters of the latter construction, one can get functions with intermediate behavior, i.e., for which the Kolmogorov--Sinai entropy vanishes, the lower permutation entropy is finite but positive, and the upper permutation entropy is infinite.

\section{Preliminaries}
\subsection{Shannon entropy}
We briefly recall the definition of Shannon entropy and a fundamental inequality regarding the entropy of mixed distributions, which will be essential for our later estimates.
\begin{definition}
    The Shannon entropy or just entropy of a discrete random variable $X$ is 
    \[
    H[X] = - \sum_x\mathbb{P}(X=x)\log\mathbb{P}(X=x)=- \mathbb{E}[\log\mathbb{P}(X)].
    \]

    By the Shannon entropy of a discrete probability measure $\mu$, we mean the Shannon entropy of a random variable $X$ with $X\sim \mu$.
\end{definition}
Note that if $X$ is a discrete random variable taking values in a set of size $n$, then $H[X] \le \log n$, with equality if and only if $X$ is uniformly distributed.

We establish the following fundamental properties of Shannon entropy to be used in our entropy rate estimations:
\begin{lemma}[Refinement] \label{lem:refinement}
Let $\mathcal{P} = (p_1, \dots, p_k)$ be a probability distribution. Let $\mathcal{Q}$ be a refinement of $\mathcal{P}$ given by $\mathcal{Q} = (q_{1,1}, \dots, q_{1,l_1}; \dots; q_{k,1}, \dots, q_{k,l_k})$ such that $\sum_{j=1}^{l_i} q_{i,j} = p_i$ for each $i \in \{1, \dots, k\}$. Then:
\[
H(\mathcal{Q}) \ge H(\mathcal{P}).
\]
\end{lemma}
\begin{lemma}[Mixing lemma] \label{2}
Let $(\nu_k)_{k=1}^m$ be probability measures on the same finite ground set $E$ and $(\alpha_k)_{k=1}^m$ a nonnegative vector with $\ell_1$ norm $\|\alpha\|_1>0$. Then
\[
\sum_{k=1}^m \alpha_k H(\nu_k) \le \|\alpha\|_1 H \left( \frac{\sum_{k=1}^m \alpha_k \nu_k}{\|\alpha\|_1} \right) \le \sum_{k=1}^m \alpha_k H(\nu_k) + \|\alpha\|_1 H \left( \frac{(\alpha_k)_{k=1}^m}{\|\alpha\|_1} \right).
\]
\end{lemma}

\vspace{9pt}
\subsection{Permutation entropy for interval maps}
We now consider a classical definition of permutation entropy for interval maps. While this framework does not require the map to be measure-preserving, we shall restrict our discussion to the measure-preserving setting to maintain consistency with our subsequent constructions.

Let $f : I \to I$ be an interval map preserving the Lebesgue measure $\lambda$. For any initial point $x$, we consider its first $n$ iterates $(x, f(x), \dots, f^{n-1}(x))$. Given a permutation $\pi = (k_1, \dots, k_n)$ of $\{0, \dots, n-1\}$, we define $P_\pi$ as the set of points $x$ whose orbit is ordered according to $\pi$:$$ P_\pi = \{x \in I \mid f^{k_1}(x) < f^{k_2}(x) < \dots < f^{k_n}(x)\}.$$
For $n\geq 2$ let $\mathcal{P}_n^*$ be the collection of the $P_\pi$ which are nonempty. Note that the $P_\pi$ are disjoint and cover $[0, 1]$ apart from points for which $f^i(x)=f^{j}(x)$ for some $0\leq i<j\leq n-1$, i.e., periodic and eventually periodic points recurring in $n-1$ steps. Assume temporarily that periodic points form a null set, then $\mathcal{P}_n^*$ partitions $[0, 1]$ apart from a null set for every $n$. We define the order $n$ permutation entropy of the map $f$ as the Shannon entropy of the distribution induced by the partition $\mathcal{P}^*_n$: \begin{align*}
    \Hpn(f) &= -  \sum_{P_{\pi} \in \mathcal{P}^*_n} \lambda(P_\pi) \log (\lambda (P_{\pi})) .
\end{align*} 

We denote the probability distribution over $S_n$ induced by the partition $\mathcal{P}^*_n$ (with respect to Lebesgue measure) by $\mu_f^n$. Using this notation, we have $\Hpn(f)=H(\mu_f^n)$.
The permutation entropy of $f$ is defined as $\Hp(f)=\lim_{n\to \infty} \Hpn(f)/n$, whenever the limit exists.

The literature of permutation entropy mostly concerns with the case discussed above, where periodic points form a null set. To handle functions that may possess periodic points, we adopt the following convention: we extend the state space of ordinal patterns from the symmetric group $S_n$ to $S_n \cup \{\ast\}$, where the extra symbol $\ast$ represents the event that the sequence of iterates contains at least one equality (i.e. $f^i(x) = f^j(x)$ for some $0 \le i < j \le n-1$). We note that introducing this extra state increases the theoretical maximum of $\Hpn(f)$ to $\log(n!+1)$, but that makes no difference asymptotically. We further note that this convention is entirely harmless to our core investigation. A theorem by Bobok and Troubetzkoy \cite{BCOT22} establishes that for a generic continuous measure-preserving function in $\mathcal{C}(\lambda)$, the set of periodic points has Lebesgue measure zero. Consequently, for generic maps, the probability mass assigned to the state $\ast$ is exactly zero, and this convention does not affect the asymptotic behavior of the entropy. However, this definition allows the permutation entropy $\Hpn(f)$ to be well-defined for all $f \in \mathcal{C}(\lambda)$.

A fundamental property of permutation entropy is that it recovers the Kolmogorov-Sinai entropy of the map as $n$ grows. Specifically, for a wide class of maps, the following result holds \cite{BandtKellerPompe}:

\begin{theorem}\label{permentropythm} For any piecewise monotone interval map $f: [0,1] \to [0,1]$ with finitely many pieces which preserves $\lambda$, we have$$ \lim_{n \to \infty} \frac{1}{n} \Hpn(f) = h_{KS}(f). $$\end{theorem}

As a direct consequence of this theorem combined with Rokhlin's formula \cite{5}, we can explicitly determine the permutation entropy rate when the derivative is uniformly bounded away from one.

\begin{corollary}\label{Rokhlin} 
    Suppose $f$ is a $\lambda$-preserving piecewise monotone interval map $f: [0,1] \to [0,1]$ with finitely many pieces. If there exists $K>1$ such that $|f'| >K$ almost everywhere, the permutation entropy is given by the integral of the log-derivative:
    \begin{equation}\label{eq:rokhlin}
        \lim_{n \to \infty} \frac{1}{n} \Hpn(f) = \int_0^1 \log |f'| \, d\lambda.
    \end{equation}
\end{corollary}

Another immediate corollary of Theorem \ref{permentropythm} establishes the permutation entropy rate for irrational rotations, a fact we will require for later constructions, taking into account the well-known fact that any rotation has Kolmogorov--Sinai entropy 0 (see e.g. \cite[Proposition~4.7]{SarigNotes}).

\begin{corollary}\label{3}
    Let $f: [0,1]\to [0,1]$ be an irrational rotation, i.e., $f(x)=x+\alpha \pmod{1}$ with $\alpha \in \mathbb{R} \setminus \mathbb{Q}$. Then $\lim_{n\to\infty}\frac{1}{n} \Hpn(f) = 0.$
\end{corollary}

\vspace{9pt}
\subsection{An auxiliary lemma}
The following lemma is a standard result in dynamical systems; we provide a brief proof to make this work self-contained.

\begin{lemma} \label{uniform_iterates}
Let $I = [0,1]$ and let $(f_n)_{n \in \mathbb{N}}$ be a sequence of functions $f_n: I \to I$ that converges uniformly to a continuous function $f: I \to I$. Then for any integer $k \ge 1$, the sequence of iterates $f_n^k$ converges uniformly to $f^k$ on $I$.
\end{lemma}

\begin{proof}
We proceed by induction on $k$. The base case $k=1$ is exactly the hypothesis of uniform convergence. Assume the statement holds for $k-1$, meaning $f_n^{k-1} \to f^{k-1}$ uniformly on $I$. 

Let $\varepsilon > 0$. We need to show there exists an $N$ such that for all $n \ge N$ and all $x \in I$, $|f_n^k(x) - f^k(x)| < \varepsilon$. By the triangle inequality,
\[
|f_n^k(x) - f^k(x)| \le \left|f_n(f_n^{k-1}(x)) - f(f_n^{k-1}(x))\right| + \left|f(f_n^{k-1}(x)) - f(f^{k-1}(x))\right|.
\]
Since $I$ is compact, the continuous function $f$ is uniformly continuous on $I$. Thus, there exists a $\delta > 0$ such that for any $y, z \in I$ with $|y - z| < \delta$, we have $|f(y) - f(z)| < \varepsilon/2$. By the inductive hypothesis, there exists an $N_1$ such that for all $n \ge N_1$ and all $x \in I$, $|f_n^{k-1}(x) - f^{k-1}(x)| < \delta$. Therefore, for $n\ge N_1$, the second term in the triangle inequality is strictly less than $\varepsilon / 2.$

By the uniform convergence of $f_n \to f$, there exists an $N_2$ such that for all $n \ge N_2$ and all $y \in I$, $|f_n(y) - f(y)| < \varepsilon/2$. 

Let $N = \max(N_1, N_2)$. For all $n \ge N$ and $x \in I$, $|f_n^k(x) - f^k(x)| < \varepsilon$, completing the proof.
\end{proof}

\vspace{8pt}
\section{Explicit examples}\label{sec:examples}
\subsection{A function with infinite permutation entropy}

\begin{proof}[Proof of Proposition \ref{prop:1}] To begin, we fix $n$ and construct a function $f_n$ that maximizes the order $n$ permutation entropy. Given that $\Hpn=H(\mathcal{P}^*_n)\le \log (\lvert \mathcal{P}^*_n \rvert ) \le \log (n!)$ with equality if and only if $\mathcal{P}^*_n$ consists of $n!$ sets of equal measure, the entropy is maximized when each permutation $\pi \in S_n$ occurs with uniform probability.  

To proceed, we define the cyclic variants of a permutation $\pi = (k_1, \dots, k_n) \in S_n$ as the set$$ C(\pi) = \{ (k_1, \dots, k_n), (k_2, \dots, k_n, k_1), \dots, (k_n, k_1, \dots, k_{n-1}) \}. $$ Since each such set contains exactly $n$ distinct permutations and $C(\pi) \neq C(\pi')$ implies $C(\pi) \cap C(\pi')=\emptyset$, $S_n$ can be decomposed into $(n-1)!$ disjoint cyclic classes:$$ S_n = \bigsqcup_{i=1}^{(n-1)!} C(\pi_i) $$for a properly chosen set of representative permutations $\{\pi_1, \dots, \pi_{(n-1)!}\}$. 

We now choose a set of representative permutations $\{\pi_1, \dots, \pi_{(n-1)!}\}$ and construct the function $f_{C(\pi_i)}$ for each $i=1,\ldots, (n-1)!$ as follows. (To ease notation, we think of elements of $S_n$ as permutations of $\{0, 1, \dots, n-1\}$.)
For a given class $C(\pi_i)$, let $(k_0, k_1, \dots, k_{n-1}) \in C(\pi_i)$ be a representative starting with $k_0=0$. We define $f_{C(\pi_i)}$ such that its graph is contained by the union of $n$ subsquares $\left[\frac{k_j}{n}, \frac{k_j+1}{n}\right]\times \left[\frac{k_{j+1}}{n}, \frac{k_{j+1}+1}{n}\right], \ j=0,\ldots, n-1,$ where we define $k_n=k_0=0$; see Figure\ref{fig:combined_permuton}.

To ensure $\Hpn$ is well-defined, we construct $f_{C(\pi_i)}$ by placing identical, scaled copies of a base measure-preserving function $\tilde{f}_n$ into each of these $n$ squares. Specifically, within each chosen subsquare, the graph of $f_{C(\pi_i)}$ is a copy of $\tilde{f}_n$ that has been scaled by a factor of $1/n$ and shifted to the appropriate grid position. We will define our exact choice of $\tilde{f}_n$ later. Let  $f_{C(\pi_i)}$ be the resulting measure-preserving function on the unit interval. By this construction, if a point $x$ belongs to the interval $[\frac{k_j}{n}, \frac{k_j+1}{n}]$, its iterates satisfy $f_{C(\pi_i)}^l(x)\in [\frac{k_{j+l}}{n}, \frac{k_{j+l}+1}{n}]$. Consequently, we obtain $x\in P_{\pi^{-1}}$ with $\pi=(k_j, k_{j+1},\ldots, k_{j-1})$. Thus, choosing $x$ uniformly randomly, we obtain each $\pi^{-1}$ with $\pi \in C(\pi_i)$ with equal probability.

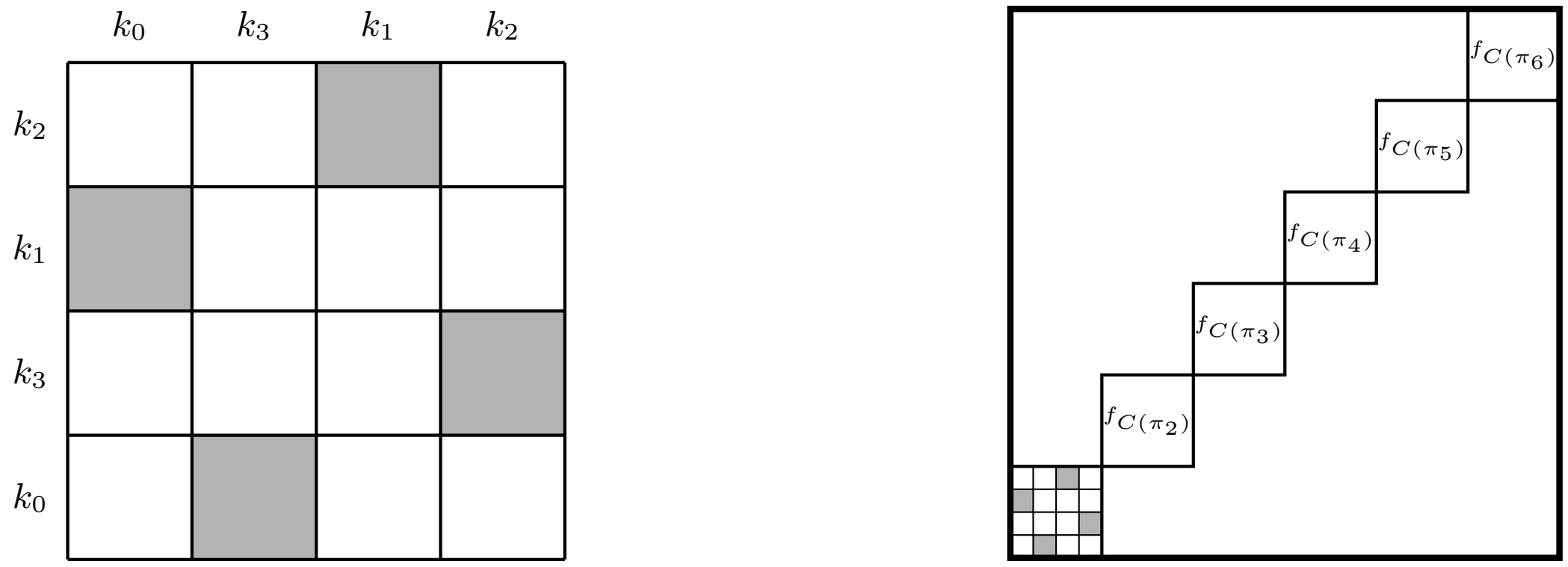
\begin{figure}[H]
    \centering
    \begin{subfigure}[b]{0.45\textwidth}
        \centering
        \begin{tikzpicture}[scale=1.1]
            \fill[gray!60] (0,2) rectangle (1,3); 
            \fill[gray!60] (1,0) rectangle (2,1); 
            \fill[gray!60] (2,3) rectangle (3,4); 
            \fill[gray!60] (3,1) rectangle (4,2); 

            \draw[thick] (0,0) grid (4,4);

            \node at (0.5, 4.3) {\small $k_0$};
            \node at (1.5, 4.3) {\small $k_3$};
            \node at (2.5, 4.3) {\small $k_1$};
            \node at (3.5, 4.3) {\small $k_2$};

            \node at (-0.3, 0.5) {\small $k_0$};
            \node at (-0.3, 1.5) {\small $k_3$};
            \node at (-0.3, 2.5) {\small $k_1$};
            \node at (-0.3, 3.5) {\small $k_2$};
        \end{tikzpicture}
    \end{subfigure}
    \hfill
    \begin{subfigure}[b]{0.45\textwidth}
        \centering
        \begin{tikzpicture}[scale=0.81] 
            \begin{scope}[shift={(0,0)}, scale=0.25]
                \fill[gray!60] (0,2) rectangle (1,3);
                \fill[gray!60] (1,0) rectangle (2,1);
                \fill[gray!60] (2,3) rectangle (3,4);
                \fill[gray!60] (3,1) rectangle (4,2);
                \draw[thin] (0,0) grid (4,4);
            \end{scope}

            \foreach \i in {0,1,2,3,4,5} {
                \draw[thick] (\i,\i) rectangle (\i+1,\i+1);
            }

            \node at (1.5,1.5) {\tiny $f_{C(\pi_2)}$};
            \node at (2.5,2.5) {\tiny $f_{C(\pi_3)}$};
            \node at (3.5,3.5) {\tiny $f_{C(\pi_4)}$};
            \node at (4.5,4.5) {\tiny $f_{C(\pi_5)}$};
            \node at (5.5,5.5) {\tiny $f_{C(\pi_6)}$};

            \draw[ultra thick] (0,0) rectangle (6,6);
        \end{tikzpicture}
    \end{subfigure}
    
    \vspace{0.3cm}
    \caption{Left: Grid representation of the function $f_{C((2,3,1,0))}$. Right: The function $f_4.$}
    \label{fig:combined_permuton}
\end{figure}

Finally, by arranging the $(n-1)!$ component functions $f_{C(\pi_i)}$ in a diagonal structure (see Figure \ref{fig:combined_permuton}) and assigning each an interval of length $\frac{1}{(n-1)!}$, we obtain the function $f_n$, where each permutation $\pi\in S_n$ is obtained with uniform probability, as required.

We are now ready to construct the function $f$. Let $(c_k)_{k=1}^{\infty}$ be a non-decreasing sequence of positive integers. We partition the unit interval $[0,1]$ by placing a sequence of intervals $I_k$ of length $2^{-k}$, starting from $0$ for $k=1$. On each interval $I_k$, we embed the function $f_{c_k}$ (scaled and shifted to map $I_k$ onto itself).
To ensure the lack of periodic points, we define the base function $\tilde{f}_{c_k}(x) = x+\alpha_{c_k} \pmod 1$ with $\alpha_{c_k}\in \mathbb{R}\setminus \mathbb{Q}$ arbitrary.

It is clear then that $h_{KS}(f)=0$: each $I_k$ is an invariant set over which $f'=1$. Thus restricted to $I_k$, even the permutation entropy is zero by Theorem \ref{permentropythm}.

For each $n$, we define the probability distributions $\mu_1^n$ and $\mu_2^n$ as the normalized sums of the components $\mu_{c_k}^n$ for indices $c_k \ge n$ and $c_k < n$, respectively, that is $\mu_1^n=a_n\cdot\sum_{k: c_k\ge n} \frac{1}{2^k} \mu_{c_k}^n$ and $\mu_2^n=b_n\cdot \sum_{k: c_k< n} \frac{1}{2^k} \mu_{c_k}^n$. This way we get:

$$ \mu^n = \left( \sum_{k: c_k\ge n} \frac{1}{2^k} \right) \mu^n_1 + \left( \sum_{k: c_k < n} \frac{1}{2^k} \right) \mu^n_2 $$
Using Lemma \ref{2}:
\begin{equation}\label{1}
    \frac{\Hp(\mu^n)}{n}\ge  \sum_{k: c_k\ge n} \frac{1}{2^k} \cdot \frac{\Hp(\mu ^n_1)}{n}+  \sum_{k: c_k< n} \frac{1}{2^k} \cdot \frac{\Hp(\mu^n_2)}{n}\ge  \sum_{k: c_k\ge n} \frac{1}{2^k} \cdot \frac{\log(n!)}{n} + 0 
\end{equation}

Here we used the fact that for each $c_k\ge n$, we obtain all permutations of $S_n$ from the component $f_{c_k}$ with the same probability.
Let us choose \[c_k={2^{2^{2^k}}},\] then $c_k\ge n \iff k\ge \log_2\log _2\log_2(n).$ Then the following lower bound holds:
\[
 \sum_{k: c_k\ge n} \frac{1}{2^k}\ge 2^{-\log_2\log _2\log_2(n)}=\frac{1}{\log _2\log_2(n)}.
\]
Substituting back to (\ref{1}) we get:
\begin{equation}
    \frac{H(\mu^n)}{n}\ge \frac{\log(n!)}{n\cdot \log_2 \log_2(n)} \ge \frac{\frac{n}{2}\log(n)}{n\cdot \log_2 \log_2(n)}=\frac{1}{2}\frac{\log(n)}{\log_2 \log_2(n)} \to \infty, \text{ as } n\to \infty.
\end{equation}

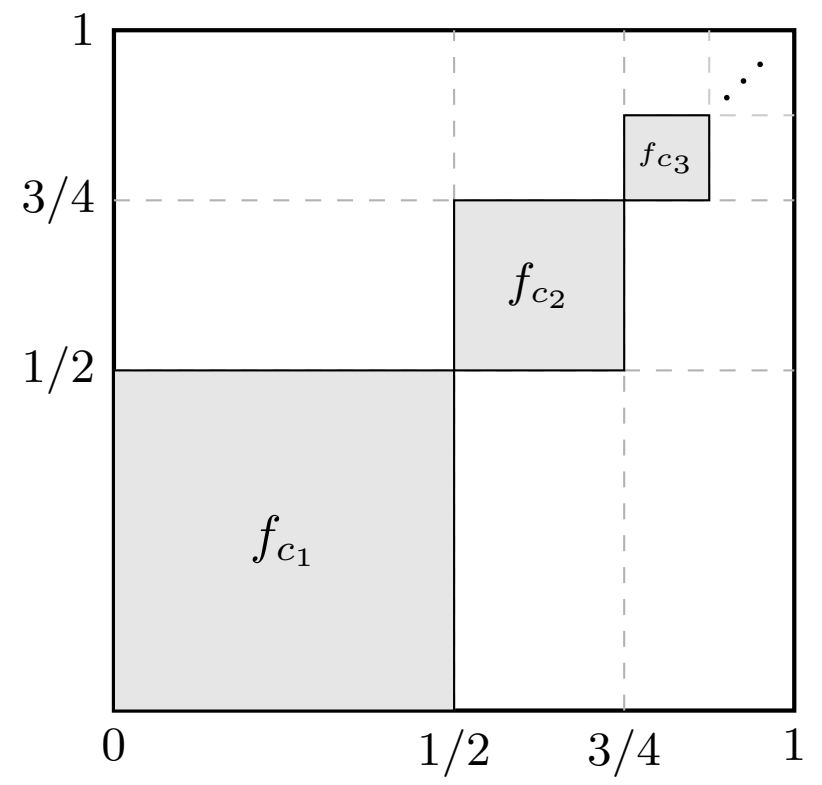
\begin{figure}[H]
    \centering
    \begin{tikzpicture}[scale=4.5] 
        \draw[thick] (0,0) rectangle (1,1);
        
        \node[below] at (0,0) {\small $0$};
        \node[below] at (0.5,0) {\small $1/2$};
        \node[below] at (0.75,0) {\small $3/4$};
        \node[below] at (1,0) {\small $1$};
        
        \node[left] at (0,0.5) {\small $1/2$};
        \node[left] at (0,0.75) {\small $3/4$};
        \node[left] at (0,1) {\small $1$};

        \draw[dashed, gray!60] (0,0.5) -- (1,0.5);
        \draw[dashed, gray!60] (0.5,0) -- (0.5,1);
        \draw[dashed, gray!60] (0,0.75) -- (1,0.75); 
        \draw[dashed, gray!60] (0.75,0) -- (0.75,1); 
        
        \draw[dashed, gray!40] (0.75,0.875) -- (1,0.875);
        \draw[dashed, gray!40] (0.875,0.75) -- (0.875,1);

        \fill[gray!20, draw=black] (0,0) rectangle (0.5,0.5);
        \node at (0.25,0.25) {$f_{c_1}$};

        \fill[gray!20, draw=black] (0.5,0.5) rectangle (0.75,0.75);
        \node at (0.625,0.625) {\small $f_{c_2}$};

        \fill[gray!20, draw=black] (0.75,0.75) rectangle (0.875,0.875);
        \node at (0.8125,0.8125) {\tiny $f_{c_3}$};

        \node[rotate=45] at (0.93, 0.93) {$\dots$};
    \end{tikzpicture}
    \caption{The dyadic-diagonal construction of the function $f$. The graphs of the components $f_{c_k}$ are embedded in squares along the diagonal corresponding to subintervals of length $2^{-k}$.}
    \label{fig:diadic_construction}
\end{figure}

\end{proof}

\vspace{8pt}
\subsection{A function with oscillating entropy rate}
\begin{proof}[Proof of Proposition \ref{prop:2}] Let us define the measure-preserving function $g_{\beta}(x)=x+\beta \pmod 1$. For each $k\in \mathbb{N}$, choose irrational $0<\beta_k\le \frac{1}{k}$, so that the first $k$ iterates of $g_{\beta_k}$ do not "overlap" for any starting point $x$. More precisely, if we fix $k$, for any $n\le k$ we have $\{\pi \in S_n: P_{\pi} \in \mathcal{P}^*_n\} = C(id),$ therefore $\lvert \mathcal{P}^*_n \rvert = n$ for all $n\le k$.

Let $\tilde{f}_{c_k}(x) = g_{\beta_k} = x + \beta_k \pmod 1 $ be the local base map used to construct the component $f_{c_k}$. This choice defines the micro-structure at level $c_k$. For $n \ge c_k$, a point's orbit follows the macro-pattern of length $c_k$ determined by the chosen cyclic permutation $\pi \in S_{c_k}$. However, the precise ordering of iterates is further refined by the local behavior of $f_{c_k}$ inside each subsquare.
To understand the structure of $\mu_{c_k}$, consider the orbit of a point under this construction. Globally, the point 'jumps' between subsquares according to the cyclic permutation $\pi \in S_{c_k}$ at each time step. Locally, however, the point’s relative position within each subsquare is governed by the map $f_{c_k}$. Because each subsquare contains an identical, scaled copy of $f_{c_k}$, the point effectively 'climbs' the graph of the function while being shifted across the grid. Specifically, let $\mathcal{P}$ be a partition formed by $c_k!$ scaled copies of the partition associated with $f_{c_k}$. For any point starting in an element $I \in \mathcal{P}$, the sequence of subsquares it visits is fixed by the cycle, and its relative vertical rank within those squares is fixed by the local dynamics of $f_{c_k}$. Thus, the resulting order-$n$ permutation is uniquely determined for all $x \in I$. This shows that $\mathcal{P}$ is a refinement of $\mathcal{P}_n^*(f_{c_k})$, allowing us to bound the entropy: $H(\mu_{c_k}) \le H(\mathcal{P}).$
By applying the mixing inequality for Shannon entropy (Lemma \ref{2}), we can bound the total entropy of the order-$n$ permutation distribution: $\Hpn(f_{c_k}) \le H(\mathcal{P})\le \Hpn(\tilde{f}_{c_k}) + \log (c_k!)$. Thus, using Corollary \ref{3}, $\frac{\Hpn (f_{c_k})}{n}\to 0,$ as $n\to \infty.$

\begin{figure}[H]
    \includegraphics[width=.75\textwidth,center]{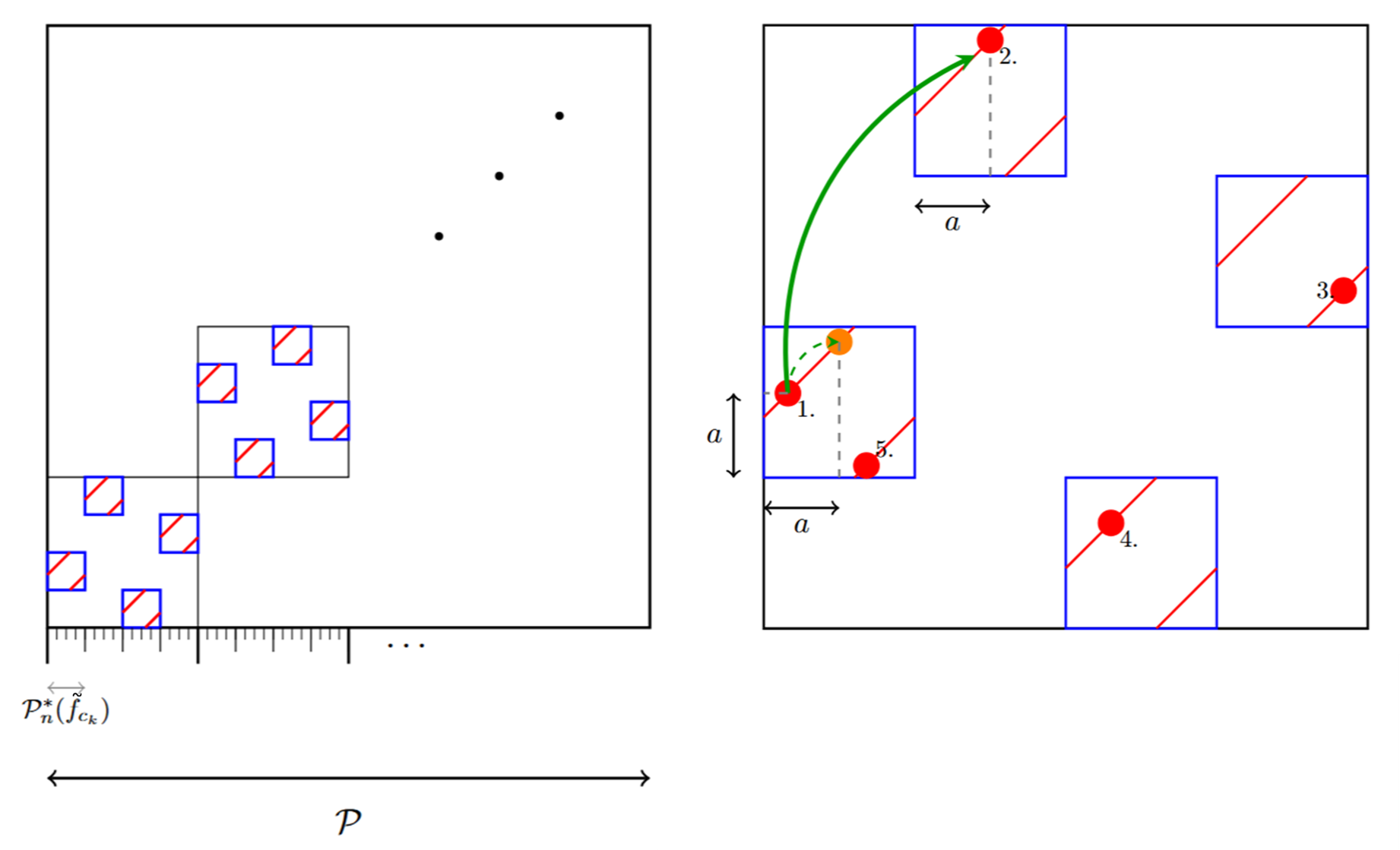}
    \caption{Left: The function $f_4$ with the rotation map embedded and a schematic representation of the partition $\mathcal{P}$; Right: First 5 iterates in the subsquare associated to $C(1,3,2,0)$.}
\end{figure}

\subsubsection*{Recursive Construction of the Function}
Finally, we are ready to construct our function $f$ with oscillating entropy behavior. We define $f$ using the same dyadic-diagonal construction described previously. To each $j$-th subinterval along the diagonal, we assign a function component scaled to a length of $2^{-j}$. The specific assignment is determined by two sequences: a sequence of indices $n_k \to \infty$ and a sequence of  indices $m_k \to \infty$, chosen such that $n_{k-1} < m_k < n_k$.

For each $k \in \mathbb{N}$, the assignments are as follows:
\begin{itemize}
\item To the index $j = n_k$, we assign the high-entropy component $f_{c_{n_k}}$.
\item To all other indices $j$ in the interval $m_k \le j < m_{k+1}$, we assign the rotational map $g_{\beta_{m_k}}$.
\end{itemize}

Observe first that regardless of the choice of $n_k, m_k$, the high-entropy components already guarantee $\limsup_{n\to\infty} \Hpn(f)/n=\infty$. Indeed, by the computations of the previous construction we can infer using $k$th high entropy part that
$$\frac{H_{c_{n_k}}^{\mathrm{perm}}(f)}{c_{n_k}}\geq 2^{-n_k}\frac{\log ({c_{n_k}}!)}{c_{n_k}}\sim 2^{-n_k}\log c_{n_k} \to \infty.$$

Consequently, it suffices to show that by carefully choosing the growth rates of $n_k$ and $m_k$, the resulting entropy rate along the sequence $m_k$ satisfies:
\[ \lim_{k \to \infty} \frac{\Hp_{m_k}(f)}{m_k} = 0. \]

\subsubsection*{Estimating the Entropy of the Preceding Indices}
We define our sequences inductively. For $k > 1$, assume that $n_{k'}$ and $m_{k'}$ have been defined for all $k' < k$. For $k=1$, we simply assume $m_1=1$ and $n_1=2$. Suppose we have fixed $m_{k}$. We seek an upper bound for the entropy contributed by the indices $j < m_{k}$. To this end, let $\gamma_{k}$ be a normalizing factor such that the weighted sum of the induced distributions forms a probability distribution $\nu_{k}^{m_{k}}$:$\nu_{k}^{m_{k}}$:$$\nu_{k}^{m_{k}} = \gamma_{k} \left[ \sum_{i=1}^{k-1} \left( \sum_{\substack{m_i \le j < m_{i+1} \\ j \neq n_i}} 2^{-j} \right) \mu_{\beta_{m_i}}^{m_{k}} + \sum_{i=1}^{k-1} 2^{-n_i} \mu_{c_{n_i}}^{m_{k}} \right].$$Applying the mixing inequality for Shannon entropy and noting that the number of distinct function components is at most $2k$, we obtain:$$H(\nu_{k}^{m_{k}}) \le \log(2k) + \sum_{i=1}^{k-1} H^{\text{perm}}_{m_{k}}(g_{\beta_{m_i}}) + \sum_{i=1}^{k-1} H^{\text{perm}}_{m_{k}}(f_{c_{n_i}}).$$Since $k$ is fixed, each term in the summands represents the permutation entropy of a fixed function component. Because these entropies vanish upon normalization by $m_k$ as the permutation order $m_{k} \to \infty$, we can choose $N$ sufficiently large such that $H(\nu_{k}^{m_{k}})/m_k < \frac{1}{3k}$ for all $m_{k} > N$.

\subsubsection*{Estimating the Entropy of the Tail}

For the indices $j \ge m_{k}$, we construct the probability measure $\rho_{k}^{m_{k}, n_{k}}$ in a similar fashion:
\[
\rho_{k}^{m_{k}, n_{k}} = \delta_{m_{k}, n_{k}} \left( a_{m_{k}, n_{k}} \cdot B_1(m_{k}, n_{k}) + b_{m_{k}, n_{k}} \cdot B_2(m_{k}, n_{k}) \right)
\]
where:
\begin{itemize}
    \item $\delta_{m_{k}, n_{k}}$ is the normalization factor.
    \item $B_1$ represents the component of the function associated with future indices $j = n_{k'}$ for all $k \le k'$.
    \item $B_2$ represents the component associated with the remaining indices $j \ge m_{k}$ (where $j \neq n_{k'}$).
    \item $a_{m_{k}, n_{k}}$ and $b_{m_{k}, n_{k}}$ are the respective total weights of these components within $f$.
\end{itemize}

Note that while the specific values of $n_{k'}$ (for $k' \ge k$) are not yet determined, we can nonetheless establish uniform upper bounds for the entropy. These bounds will depend exclusively on the value of $n_{k}$ and $m_{k}$. 

To establish an upper bound, we recall that the rotation parameters $\beta_{k'}$ are chosen such that the permutation partition satisfies $\mathcal{P}^*_{m_{k}}(B_2) = C(\text{id})$. This ensures that the iterates within the $B_2$ component exhibit minimal combinatorial complexity. For the $B_1$ component, we choose $n_{k}$ sufficiently large so that the weight $a_{m_{k}, n_{k}}$ effectively suppresses the high-entropy terms. Specifically, we have:$$\frac{H^{\perm}_{m_{k}} (\rho_{k}^{m_{k}, n_{k}})}{m_{k}} \le \frac{\log 2}{m_{k}} + \frac{\log(m_{k}!)}{m_{k} \cdot 2^{n_{k}-m_k}} + \frac{\log m_{k}}{m_{k}} \le \frac{\log 2}{m_{k}} + \frac{\log(m_{k})}{2^{n_{k}-m_k}} + \frac{\log m_{k}}{m_{k}}.$$The first and third terms depend only on $m_{k}$; thus, we can choose $N < m_{k}$ such that their sum is less than $\frac{1}{3k}$. Finally, by selecting $n_{k} > m_{k}$ sufficiently large, the middle term also falls below $\frac{1}{3k}$, yielding the desired bound.

\vspace{7 pt}
Combining the estimates we obtain the global bound for the order-$m_{k}$ permutation entropy of $f$ \[
\frac{\Hp_{m_{k}} (f)}{m_{k}}\le \frac{1}{k}+\frac{\log 2}{m_{k}} \to 0, \text{ as } k\to \infty.
\]

This confirms that the function $f$ exhibits oscillating entropy behavior, where the entropy rate vanishes along the sequence $m_{k}$ but diverges along $c_{n_{k}}$.

\end{proof}

\section{Generic Behavior of Permutation Entropy}\label{sec:main_result_proof}

In this section, we present the proof of our main result, Theorem \ref{mainthm}. We begin by establishing a crucial lemma regarding the continuity of the permutation entropy, which will allow us to deduce the topological openness of the relevant subsets of functions.

\begin{lemma} \label{lemma3}
Let \((f_n)_{n \in \mathbb{N}}\) be a sequence of continuous measure-preserving maps on \(I = [0,1]\) converging uniformly to the continuous measure preserving map \(f\).  Suppose that for a given integer $k \ge 2$, the set of points with period strictly less than $k$ has Lebesgue measure zero for the limit function $f$, that is:
\[
\lambda\left( \bigcup_{i=1}^{k-1} \{x \in I \mid f^i(x) = x\} \right) = 0.
\]
Then, the order $k$ permutation entropy functional is continuous at $f$:
\[
\lim_{n \to \infty} \Hp_k(f_n) = \Hp_k(f).
\]
\end{lemma}

\begin{proof}
Let \(F_f = \bigcup_{0 \le i < j \le k-1} \{x \in I \mid f^i(x) = f^j(x)\}\) denote the set of points where the first \(k\) iterates of \(f\) are not pairwise distinct. Because \(F_f = \bigcup_{0 \le i < j \le k-1} (f^i)^{-1}(\text{Fix}(f^{j-i}))\), and \(f\) is measure-preserving with \(\lambda(\text{Fix}(f^m)) = 0\) for all \(1 \le m < k\) by hypothesis, it follows that \(\lambda(F_f) = 0\). From the continuity of \(f\), we get that \(F_f\) is closed, hence its complement \(G_f = I \setminus F_f\) is an open set of full measure, \(\lambda(G_f) = 1\).

Since the Shannon entropy function is continuous on the finite probability simplex, establishing convergence of the distribution \(\mu_f^k\) in the \(L^1\) norm, \(\lim_{n \to \infty} \|\mu_{f_n}^k - \mu_f^k\|_1 = 0\), is sufficient to guarantee the convergence of the permutation entropies. Thus, we must show that for any given $\varepsilon > 0$, there exists an integer $N$ such that for all \(n \ge N\):
\[
\sum_{\pi \in S_k} \left| \lambda(P_\pi(f_n)) - \lambda(P_\pi(f)) \right| < \varepsilon.
\]

For \(x \in I\), define the minimum gap between any two distinct iterates up to order \(k-1\):
\[
\Delta^k_f(x) = \min \{|f^i(x) - f^j(x)| : i \neq j, \quad i,j \in \{0, \dots, k-1\} \}.
\]
Let \(I_m = \{x \in G_f : \Delta^k_f(x) \ge 1/m\}\). Since these sets are monotonically increasing inclusion-wise, and \(\bigcup_{m=1}^\infty I_m = G_f\), we have $\lim_{m\to \infty} \lambda \big( I_m\big)=1$. Therefore, we can choose an integer \(M\) sufficiently large such that \(\lambda(I_M) > 1 - \varepsilon/2\).

By Lemma \ref{uniform_iterates}, the uniform convergence $f_n \to f$ implies the uniform convergence of the iterates $f_n^i \to f^i$ for all $1 \le i \le k-1$. Therefore, for $\eta = \frac{1}{2M}$, there exists an integer $N$ such that for all $n \ge N$, all $x \in I$, and all $i \in \{1, \dots, k-1\}$, we have 
\[
|f_n^i(x) - f^i(x)| < \eta = \frac{1}{2M}.
\]

For any \(x \in I_M\) and \(n \ge N\), the maximum perturbation of any iterate is strictly less than half the minimum separation between the iterates of \(f\), as \(\frac{1}{2M} \le \frac{1}{2} \Delta^k_f(x)\). Consequently, the strict relative order of the sequence \((x, f_n(x), \dots, f_n^{k-1}(x))\) remains identical to the sequence \((x, f(x), \dots, f^{k-1}(x))\). Thus, for all \(x \in I_M\), \(x \in P_\pi(f_n)\) if and only if \(x \in P_\pi(f)\).

The symmetric difference \(P_\pi(f_n) \Delta P_\pi(f)\) for any permutation is therefore restricted entirely to the complement space \(I \setminus I_M\). Summing the differences across all permutations yields:
\[
\sum_{\pi \in S_k} \left| \lambda(P_\pi(f_n)) - \lambda(P_\pi(f)) \right| \le 2\lambda(I \setminus I_M) < 2\left(\frac{\varepsilon}{2}\right) = \varepsilon.
\]
Consequently, the continuity of the Shannon entropy implies that \(\lim_{n \to \infty} \Hp_k(f_n) = \Hp_k(f)\).
\end{proof}

We are now ready to prove Theorem \ref{mainthm}.
Recall that $C(\lambda)$ denotes the space of continuous measure-preserving functions on $[0,1]$, and $C^*(\lambda) \subset C(\lambda)$ is the subspace of functions whose periodic points form a set of measure zero.

The proof follows the strategy introduced in the proof of \cite[Theorem~1.4]{MB} for the analogous statement on sampling entropy. Throughout the construction, we will explicitly indicate whenever our steps rely directly on that framework, and we refer the reader to \cite{MB} for the associated technical details that carry over identically. The adaptation to permutation entropy, however, requires two major modifications. First, to ensure our approximating functions strictly belong to $C^*(\lambda)$ and satisfy the derivative bounds required to apply Rokhlin's formula (Corollary \ref{Rokhlin}), we must introduce a localized geometric modification to eliminate any affine segments with slope $\pm 1$. Second, while establishing the density of the set of functions with high entropy at certain levels was straightforward in the sampling entropy setting, it presents a substantially different challenge here, requiring an entirely distinct construction.

\begin{proof}[Proof of Theorem \ref{mainthm}]
Let $X_{N,\varepsilon} \subset C(\lambda)$ (resp. $Y_{N,\varepsilon} \subset C(\lambda)$) be the set of functions for which there exists $n > N$ with $\Hpn < \varepsilon n$ (resp. $\Hpn > (1-\varepsilon)n\log n$), that is: $$X_{N,\varepsilon} =\{f\in  C^*(\lambda): \exists n>N, \Hpn (f)<\varepsilon n\}=\bigcup_{n=N+1}^{\infty} \{f\in  C^*(\lambda): \Hpn (f)<\varepsilon n\},$$ $$Y_{N,\varepsilon} =\{f\in  C^*(\lambda): \exists n>N, \Hpn (f)>(1-\varepsilon) n \log n\}=\bigcup_{n=N+1}^{\infty} \{f\in  C^*(\lambda): \Hpn (f)>(1-\varepsilon) n \log n\}.$$ Proving that these sets are dense and open in $C(\lambda)$ would conclude the proof. Since $C^*(\lambda)$ is a residual subset of $C(\lambda)$ \cite{BCOT22}, it is sufficient to prove that $X_{N,\varepsilon}$ and $Y_{N,\varepsilon}$ are dense, relatively open subsets of $C^*(\lambda)$. Note that by Lemma \ref{lemma3}, the sets $\{f \in C^*(\lambda) : \Hpn(f) < \varepsilon n\}$ (resp. $\{f \in C^*(\lambda) : \Hpn(f) > (1-\varepsilon) n \log n\}$) are open in $C^*(\lambda)$, thus, $X_{N,\varepsilon}$ and $Y_{N,\varepsilon}$ are indeed open subsets of $C^*(\lambda)$. To establish density, we will proceed by proving that they are dense in the entire space $C(\lambda)$, which trivially implies their density in $C^*(\lambda)$.

\vspace{8pt}
We start by proving the density of $X_{N, \epsilon}\subset C(\lambda).$ Consider the dense subspace of continuous
piecewise affine measure-preserving maps, $\mathrm{CPA}(\lambda)\subset C(\lambda).$
  The idea, which we adapt from \cite{MB}, is the following: given an arbitrary $f \in \mathrm{CPA}(\lambda)$, we approximate it with a piecewise affine, measure-preserving map $g \in \mathrm{PA}(\lambda)$ being arbitrarily close to $f$. While $g$ is typically discontinuous, it is constructed to satisfy $\vert{}g'\vert{} = 1$ on all its intervals of continuity. Next, we perturb $g$ by slightly contracting the domain of each affine piece while maintaining its original range, which strictly increases the magnitude of the derivatives. This contraction naturally introduces small gaps in the domain. By bridging these gaps with properly chosen affine segments, we obtain a continuous, measure-preserving function $h_0 \in \mathrm{CPA}(\lambda)$.

To ensure our final function $h$ strictly belongs to $\mathrm{CPA}^*(\lambda)=\mathrm{CPA}(\lambda)\cap C^*(\lambda)$, we must introduce an additional modification to $h_0$ to guarantee that $\vert{}h'\vert{} \neq 1$ everywhere. If our initial construction leaves any intervals where \(|h_0'| = 1\), we modify \(h_0\) slightly further: we again shrink the domain of these specific pieces, which slightly increases their derivative. In the newly created gaps, we insert a piecewise affine ``V''-shaped segment that covers the entire vertical range of that particular piece. The resulting function \(h\) is now in \(\mathrm{CPA}^*(\lambda)\) because it is piecewise affine and its derivative never equals 1. Furthermore, since \(h\) consists of finitely many linear pieces, there exists some constant \(K > 1\) such that \(|h'| \ge K\) everywhere. This strict lower bound allows us to apply Rokhlin's formula. Intuitively, because the derivative remains very close to 1 on the vast majority of the domain and is only large within the tiny inserted gaps, Rokhlin's formula tells us that by choosing the size of these remaining gaps to be small enough, the resulting entropy \(H_n/n\) can be made arbitrarily small. Therefore, by selecting the parameters carefully throughout this construction, \(h\) successfully acts as a witness demonstrating that \(X_{N,\varepsilon}\) is dense.

We now formalize this construction. First, we must ensure the existence of a suitable function $g \in B(f,\delta)$. To do this, we partition the range $[0,1]$ into a finite collection of disjoint intervals $J_1, \ldots, J_l$, each with length strictly less than $\delta$. Because $f$ is piecewise affine, we can construct this partition such that the endpoints of the intervals $(J_j)_{j=1}^l$ include all local extrema of $f$. This ensures that the preimage $f^{-1}(J_j)$ consists of a finite number of affine subintervals, and that $f$ maps each of these subintervals surjectively onto $J_j$. Since the construction is symmetric across all $J_j$, we simplify our notation by focusing on a single interval $J = J_1$. Let $I_1, \ldots, I_k$ denote the affine subintervals forming $f^{-1}(J)$. Because $f$ is measure-preserving, we can define $g$ on each interval $I_m$ as a single affine segment with slope $\pm 1$, matching the sign of $f'$, such that their images exhaust $J$, yielding a function $g \in \mathrm{PA}(\lambda) \cap B(f,\delta)$. Next, we subdivide each domain $I_m$ into three consecutive, non-overlapping closed subintervals, denoted from left to right as $I_m^-, I_m^0,$ and $I_m^+$. For a parameter $\alpha > 0$ (which we will specify later), we uniquely define the lengths of the outer segments as:$$\vert{}I_m^-\vert{} = \alpha\vert{}I_m\vert{} \frac{\sum_{i<m}\vert{}I_i\vert{}}{\vert{}J\vert{}} \quad \text{and} \quad \vert{}I_m^+\vert{} = \alpha\vert{}I_m\vert{} \frac{\sum_{i>m}\vert{}I_i\vert{}}{\vert{}J\vert{}}.$$With the domains established, we construct the intermediate function $h_0$ to be affine on each of these three segments, preserving the local monotonicity of $f$. Specifically, on the central segment $I_m^0$, we set $h_0$ to match the image of $g$, meaning $h_0(I_m^0) = g(I_m)$. On the two outer segments, $h_0$ is defined to so that $h_0(I_m) = f(I_m) = J$. Consequently, the derivative on the outer segments satisfies $\vert{}h_0'\vert{} = \frac{\vert{}J\vert{}}{\alpha\vert{}I_m\vert{}}$. Furthermore, the length of the central segment is given by:$$\vert{}I_m^0\vert{} = \vert{}I_m\vert{}\left(1 - \alpha \frac{\sum_{i \neq m}\vert{}I_i\vert{}}{\vert{}J\vert{}}\right) = \vert{}I_m\vert{}\left(1 - \alpha\frac{\vert{}J\vert{} - \vert{}I_m\vert{}}{\vert{}J\vert{}}\right).$$ It is immediate that \(h_0\) is continuous, belongs to \(B(f,\delta)\), and is measure-preserving by the exact same algebraic verification as for \(h\) in the original construction (see the proof of Theorem 1.4 in \cite{MB}, up to this point we followed that construction).

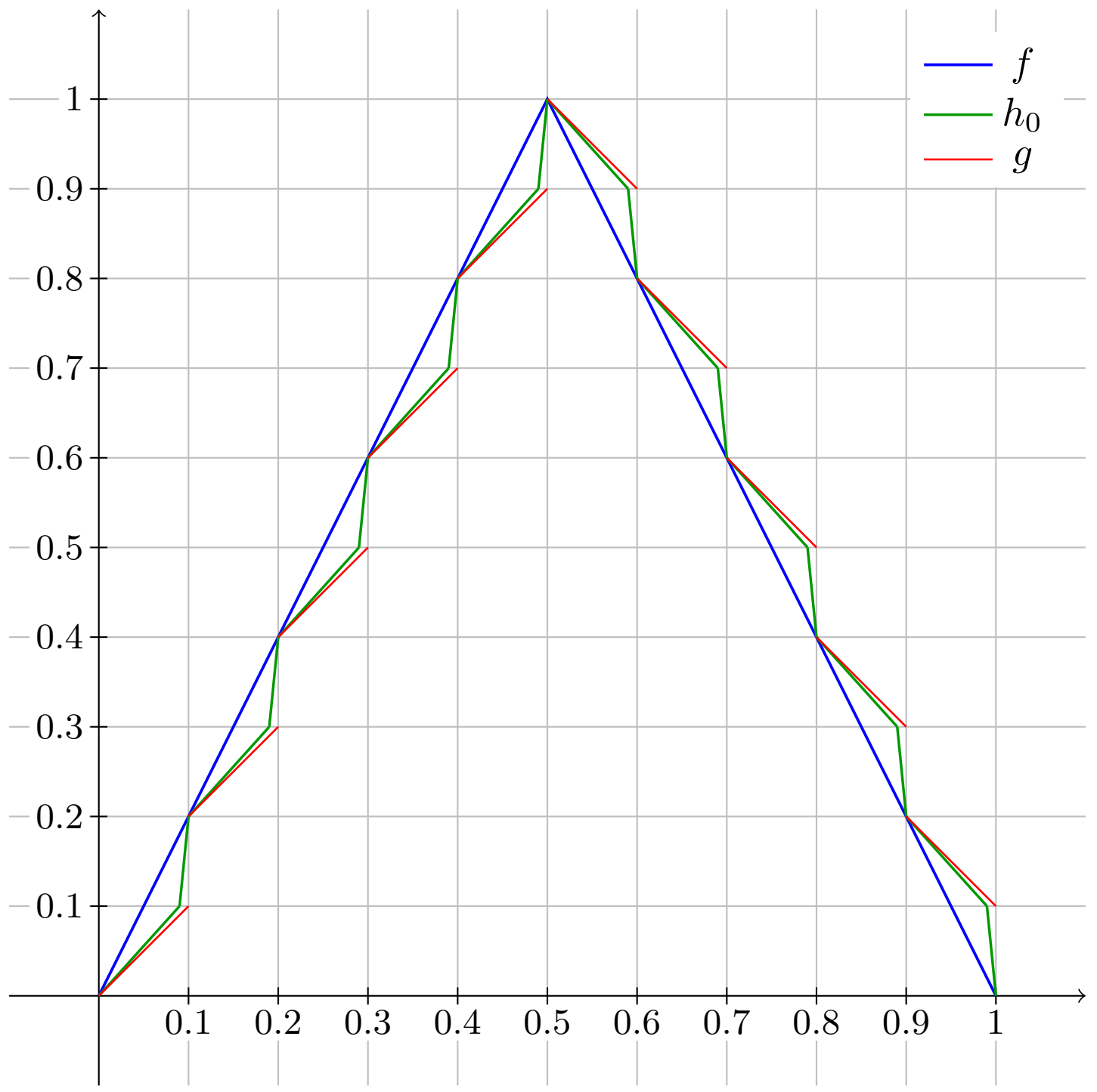
\begin{figure}[htbp]
    \centering
    \begin{tikzpicture}
\begin{axis}[
    width=11cm,
    height=11cm,
    xmin=-0.1,
    xmax=1.1,
    ymin=-0.1,
    ymax=1.1,
    axis lines=middle,
    axis line style={->},
    xtick={0,0.1,0.2,0.3,0.4,0.5,0.6,0.7,0.8,0.9, 1},
    ytick={0,0.1,0.2,0.3,0.4,0.5,0.6,0.7,0.8,0.9, 1},
    grid=both,
    tick style={
    draw=black,
    line width=0.4pt,
    preaction={
        draw=white,
        line width=2pt
    }},
    tick label style={
        font=\small,
        fill=white,
        inner sep=1.5pt
    },
    legend style={
        at={(0.98,0.98)},
        anchor=north east,
        draw=none,
        fill=white
    },
    clip=false
]

\addplot[
    blue,
    line width=0.7pt
]
coordinates {
    (0,0)
    (0.5,1)
    (1,0)
};
\addlegendentry{$f$}
\addplot[
    green!60!black,
    line width=0.65pt
]
coordinates {(0.00,0.00) (0.09,0.10)};
\addlegendentry{$h_0$}

\addplot[
    green!60!black,
    line width=0.65pt,
    forget plot
]
coordinates {(0.09,0.10) (0.10,0.20)};

\addplot[
    green!60!black,
    line width=0.65pt,
    forget plot
]
coordinates {(0.10,0.20) (0.19,0.30)};

\addplot[
    green!60!black,
    line width=0.65pt,
    forget plot
]
coordinates {(0.19,0.30) (0.20,0.40)};

\addplot[
    green!60!black,
    line width=0.65pt,
    forget plot
]
coordinates {(0.20,0.40) (0.29,0.50)};

\addplot[
    green!60!black,
    line width=0.65pt,
    forget plot
]
coordinates {(0.29,0.50) (0.30,0.60)};

\addplot[
    green!60!black,
    line width=0.65pt,
    forget plot
]
coordinates {(0.30,0.60) (0.39,0.70)};

\addplot[
    green!60!black,
    line width=0.65pt,
    forget plot
]
coordinates {(0.39,0.70) (0.40,0.80)};

\addplot[
    green!60!black,
    line width=0.65pt,
    forget plot
]
coordinates {(0.40,0.80) (0.49,0.90)};

\addplot[
    green!60!black,
    line width=0.65pt,
    forget plot
]
coordinates {(0.49,0.90) (0.50,1.00)};

\addplot[
    green!60!black,
    line width=0.65pt,
    forget plot
]
coordinates {(0.50,1.00) (0.59,0.90)};

\addplot[
    green!60!black,
    line width=0.65pt,
    forget plot
]
coordinates {(0.59,0.90) (0.60,0.80)};

\addplot[
    green!60!black,
    line width=0.65pt,
    forget plot
]
coordinates {(0.60,0.80) (0.69,0.70)};

\addplot[
    green!60!black,
    line width=0.65pt,
    forget plot
]
coordinates {(0.69,0.70) (0.70,0.60)};

\addplot[
    green!60!black,
    line width=0.65pt,
    forget plot
]
coordinates {(0.70,0.60) (0.79,0.50)};

\addplot[
    green!60!black,
    line width=0.65pt,
    forget plot
]
coordinates {(0.79,0.50) (0.80,0.40)};

\addplot[
    green!60!black,
    line width=0.65pt,
    forget plot
]
coordinates {(0.80,0.40) (0.89,0.30)};

\addplot[
    green!60!black,
    line width=0.65pt,
    forget plot
]
coordinates {(0.89,0.30) (0.90,0.20)};

\addplot[
    green!60!black,
    line width=0.65pt,
    forget plot
]
coordinates {(0.90,0.20) (0.99,0.10)};

\addplot[
    green!60!black,
    line width=0.65pt,
    forget plot
]
coordinates {(0.99,0.10) (1.00,0.00)};


\addplot[
    red,
    line width=0.5pt
]
coordinates {(0.0,0.0) (0.1,0.1)};
\addlegendentry{$g$}

\addplot[
    red,
    line width=0.5pt,
    forget plot
]
coordinates {(0.1,0.2) (0.2,0.3)};

\addplot[
    red,
    line width=0.5pt,
    forget plot
]
coordinates {(0.2,0.4) (0.3,0.5)};

\addplot[
    red,
    line width=0.5pt,
    forget plot
]
coordinates {(0.3,0.6) (0.4,0.7)};

\addplot[
    red,
    line width=0.5pt,
    forget plot
]
coordinates {(0.4,0.8) (0.5,0.9)};

\addplot[
    red,
    line width=0.5pt,
    forget plot
]
coordinates {(0.5,1.0) (0.6,0.9)};

\addplot[
    red,
    line width=0.5pt,
    forget plot
]
coordinates {(0.6,0.8) (0.7,0.7)};

\addplot[
    red,
    line width=0.5pt,
    forget plot
]
coordinates {(0.7,0.6) (0.8,0.5)};

\addplot[
    red,
    line width=0.5pt,
    forget plot
]
coordinates {(0.8,0.4) (0.9,0.3)};

\addplot[
    red,
    line width=0.5pt,
    forget plot
]
coordinates {(0.9,0.2) (1.0,0.1)};

\end{axis}
\end{tikzpicture}
    
    \caption{An illustration of the functions $g$ (red) and $h_0$ (green) constructed from the symmetric tent map $f$ (blue). The interval $J_i$ ($i = 1,\dots,5$) is fixed to be $J_i = [0.2 \cdot (i-1), 0.2 \cdot i]$, while $\alpha = 0.2$. This is Figure 4 from \cite{MB}.}
    \label{fig1}
\end{figure}

Now, we modify \(h_0\) to obtain our final function \(h \in \mathrm{CPA}^*(\lambda)\) by eliminating any pieces where \(|h_0'| = 1\). Suppose $h_0'=-1$ on some interval. By construction, for small $\alpha$, $h_0' = -1$ can only occur on a central interval $I_m^0$, and only if the original function $f$ satisfied $f' =- 1$ on $I_m$. Because $f$ is continuous and measure-preserving, any affine piece with a slope of $-1$ must lie perfectly on the anti-diagonal $y = 1 - x$.
Because \(f\) has a slope of \(-1\) on \(I_m\), it maps \(I_m\) onto the target interval \(J\) such that their lengths are exactly equal, \(|I_m| = |J|\). However, because \(f\) is measure-preserving, the sum of the lengths of all affine preimages of \(J\) must equal the length of \(J\). This immediately forces \(I_m\) to be the unique preimage of \(J\). Consequently, in our definition of the subinterval lengths, the sums representing the other preimages (\(\sum_{i<m} |I_i|\) and \(\sum_{i>m} |I_i|\)) are empty and evaluate to zero. This yields \(|I_m^-| = |I_m^+| = 0\), meaning the central interval spans the entire domain: \(I_m^0 = I_m\). 

Let \(I_m=I = [a, b]\) be such an interval with range \(J = [c, d]\).   By possibly partitioning $I$ further, we may assume $b - a < \delta$ to ensure our modifications keep $h$ within the $\delta$-neighborhood of $f$. To eliminate this derivative, we partition \(I\) into a left interval \(I_V = [a, a + \alpha|I|]\) and a right interval \(I_R = [a + \alpha|I|, b]\). We further divide \(I_V\) into \(I_{V1}\) and \(I_{V2}\), each of length \(\frac{\alpha}{2}|I|\). We define \(h\) on \(I\) as follows (see Figure \ref{fig:tent_modification_unit_square}):
\begin{itemize}
    \item On \(I_{V1}\), \(h\) decreases linearly from \(d\) to \(c\) (derivative: \(-2/\alpha\)).
    \item On \(I_{V2}\), \(h\) increases linearly from \(c\) to \(d\) (derivative: \(2/\alpha\)).
    \item On \(I_R\), \(h\) decreases linearly from \(d\) to \(c\) (derivative: \(-1/(1-\alpha)\)).
\end{itemize}

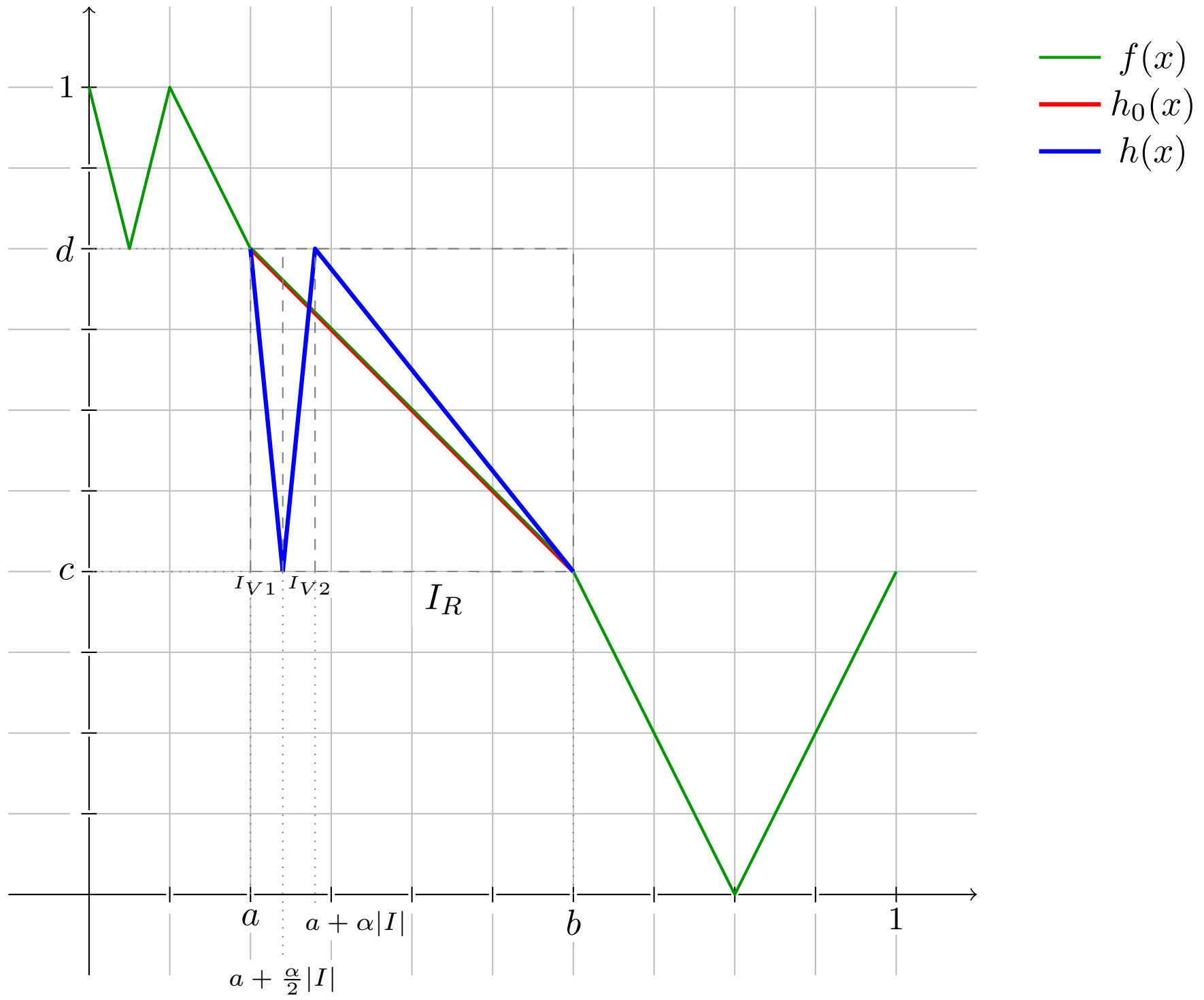
\begin{figure}[htbp]
    \centering
    \begin{tikzpicture}
    \begin{axis}[
        width=11cm,
        height=11cm,
        xmin=-0.1,
        xmax=1.1,
        ymin=-0.1,
        ymax=1.1,
        axis lines=middle,
        axis line style={->},
        xtick={0, 1},
        ytick={0, 1},
        extra x ticks={0.1,0.2,0.3,0.4,0.5,0.6,0.7,0.8,0.9},
        extra x tick labels={},
        extra y ticks={0.1,0.2,0.3,0.4,0.5,0.6,0.7,0.8,0.9},
        extra y tick labels={},
        grid=both,
        tick style={
            draw=black,
            line width=0.4pt,
            preaction={
                draw=white,
                line width=2pt
            }
        },
        tick label style={
            font=\small,
            fill=white,
            inner sep=1.5pt
        },
        legend style={
            at={(1.05, 0.98)},
            anchor=north west,
            draw=none,
            fill=white
        },
        clip=false
    ]

    \addplot[
        thick, 
        green!60!black, 
        line join=round
    ] coordinates {
        (0, 1) (0.05, 0.8) (0.1, 1) (0.2, 0.8) (0.6, 0.4) (0.8, 0) (1, 0.4)
    };
    \addlegendentry{$f(x)$}

    \addplot[
        line width=1.2pt, 
        red
    ] coordinates {
        (0.2, 0.8) (0.6, 0.4)
    };
    \addlegendentry{$h_0(x)$}

    \begin{scope}
        \clip (axis cs:0.2, 0.8) -- (axis cs:0.6, 0.4) -- (axis cs:0.6, 1.2) -- (axis cs:0.2, 1.2) -- cycle;
        \addplot[
            line width=1.2pt, 
            green!60!black, 
            forget plot
        ] coordinates {
            (0.2, 0.8) (0.6, 0.4)
        };
    \end{scope}

    \addplot[
        very thick, 
        blue, 
        line join=round
    ] coordinates {
        (0.2, 0.8) (0.24, 0.4) (0.28, 0.8) (0.6, 0.4)
    };
    \addlegendentry{$h(x)$}
    

    
    \draw[dashed, gray] (axis cs:0.2, 0.4) rectangle (axis cs:0.6, 0.8);
    
    \draw[dashed, gray] (axis cs:0.24, 0.4) -- (axis cs:0.24, 0.8);
    \draw[dashed, gray] (axis cs:0.28, 0.4) -- (axis cs:0.28, 0.8);

    \node[below left, inner sep=1pt, fill=white, text=black] at (axis cs:0.24, 0.4) {\tiny $I_{V1}$};
    \node[below right, inner sep=1pt, fill=white, text=black] at (axis cs:0.24, 0.4) {\tiny $I_{V2}$};
    \node[below, fill=white, text=black] at (axis cs:0.44, 0.4) {$I_R$};

    
    \draw[gray, dotted] (axis cs:0.2, 0.4) -- (axis cs:0.2, 0) 
        node[below, fill=white, text=black, inner sep=2pt, yshift=-2pt] {$a$};
        
    \draw[gray, dotted] (axis cs:0.6, 0.4) -- (axis cs:0.6, 0) 
        node[below, fill=white, text=black, inner sep=2pt, yshift=-2pt] {$b$};
        
    \draw[gray, dotted] (axis cs:0.24, 0.4) -- (axis cs:0.24, -0.08) 
        node[below, fill=white, text=black, inner sep=2pt] {\scriptsize $a+\frac{\alpha}{2}|I|$};
        
    \draw[gray, dotted] (axis cs:0.28, 0.4) -- (axis cs:0.28, -0.015) 
        node[below right, fill=white, text=black, inner sep=1pt, xshift=-4pt] {\scriptsize $a+\alpha|I|$};

    \draw[gray, dotted] (axis cs:0.2, 0.4) -- (axis cs:0, 0.4) 
        node[left, fill=white, text=black, inner sep=2pt, xshift=-2pt] {$c$};
        
    \draw[gray, dotted] (axis cs:0.2, 0.8) -- (axis cs:0, 0.8) 
        node[left, fill=white, text=black, inner sep=2pt, xshift=-2pt] {$d$};
    \end{axis}
    \end{tikzpicture}
    \caption{The local modification $h$ on the interval $I=[a,b]$. The overall measure-preserving map $f$ (green) and the segment $h_0 = -id$ (red) coincide perfectly on the reversed diagonal over the $4 \times 4$ grid box $I \times J$, depicted as a two-toned line. By choosing a small $\alpha$, the modified function $h$ (blue) creates steep local derivatives over the narrow intervals $I_{V1}$ and $I_{V2}$, before connecting continuously to $c$ over $I_R$.}
    \label{fig:tent_modification_unit_square}
\end{figure}

For $h_0'=1$ we make changes analogously with derivatives with opposite signs.
By replacing \(h_0\) with \(h\) on all such pieces, we ensure \(|h'| \neq 1\) everywhere. Since \(h\) consists of finitely many affine segments, each satisfying \(|h'| > 1\) (for all $\alpha>0$), there exists \(K_{\alpha}=K > 1\) such that \(|h'| \ge K\) everywhere, satisfying the prerequisite for applying Rokhlin's formula (see Corollary \ref{Rokhlin}).

To see that \(h\) remains measure-preserving, consider any subinterval \(L \subset J\) for any such modified interval \(I\). Because \(h\) is affine on each of the three sub-pieces of \(I\), the lengths of the preimages sum exactly to 
\[
|h^{-1}(L)| = |L|\left(\frac{\alpha}{2} + \frac{\alpha}{2} + (1-\alpha)\right) = |L|.
\]
Thus, modifying \(h_0\) on all such intervals preserves the measure, yielding \(h \in \mathrm{CPA}^*(\lambda)\).

Finally, we calculate the contribution of this modified piece \(I\) to the total entropy using Rokhlin's formula, \(\int \log|h'|\). The integral over \(I\) splits into the ``V''-shape and the remaining segment:
\begin{align*}
\int_I \log|h'| &= \int_{I_V} \log\left(\frac{2}{\alpha}\right) + \int_{I_R} \log\left(\frac{1}{1-\alpha}\right) \\
&= \alpha|I|\log\left(\frac{2}{\alpha}\right) + (1-\alpha)|I|\log\left(\frac{1}{1-\alpha}\right).
\end{align*}
As \(\alpha \to 0\), both terms tend to $0$. Therefore, the entropy contributed by the modified \(-1\) slope regions vanishes in the limit, and the total entropy remains bounded by \(\varepsilon n\) for sufficiently small \(\alpha\).

\vspace{10pt}
Establishing the density of $Y_{N, \epsilon}\subset C(\lambda)$ requires a more delicate construction.
We will follow an analogous construction for $h_0$ as before. We may assume without loss of generality that the affine pieces of $f$ have rational endpoints, as such functions form a dense subset of \(\mathrm{CPA}(\lambda)\) \cite{BT2020}. 

Consequently, we can fix a uniform grid on \([0,1]\) with step size \(p=\frac{1}{z}\) with $z\in \mathbb{Z^+}$, such that the partition intervals \(I\) and \(J\) align with the grid points \(\{ip\}_{i=0}^{1/p}\) in the following way.
First, we form an initial partition of the $y$-axis using the $y$-coordinates of all local extrema and endpoints of $f$. By including additional rational points if necessary, we ensure the length of each target interval $J$ is at most $\delta$. We then pull these values back to the $x$-axis via $f^{-1}$, exactly as before. Because $f$ is piecewise affine with rational coordinates, this pullback creates a partition of the domain $[0,1]$ consisting entirely of rational points, yielding subintervals $I$ on which $f$ is strictly affine. Next, we choose a sufficiently fine uniform grid on $[0,1]$ with a rational step size $p = \frac{1}{z}$ (for some $z \in \mathbb{Z}^+$). Because all boundaries of our $I$ and $J$ intervals are rational, we can choose $z$ as a common denominator to ensure that every point in our $x$-axis and $y$-axis partitions aligns with the grid points $\{ip\}_{i=0}^z$.

Using this partition structure, we first construct the intermediate function $h_0$ exactly as in the previous section with some $\alpha>0$ to be fixed later. We then construct our final function $h$ by modifying $h_0$ inside the central intervals $I_m^0$.  Our goal is to introduce $L$ rapid oscillations within each grid cell while maintaining the measure-preserving property.

First, we define the necessary structural parameters globally. Let $\mathcal{G} = \{ip : i = 0, 1, \dots, z\}$ denote the uniform grid points on the domain. Because the target intervals form a partition of the $y$-axis, each horizontal grid value $jp$ belongs to exactly one interval $J$ (up to shared boundaries). Consequently, its preimages under $h_0$ will only ever lie within the specific domain intervals $I_m$ associated with that particular $J$. This ensures our indexing is well-defined: for each $j$, if $h_0^{-1}(jp) \cap I_m^0\neq \emptyset$, let $x_{j}^m$ be its single element, i.e., the point where $h_0$ crosses the horizontal grid lines $y = jp$ (possibly coinciding with an endpoint of $I_m^0$). We then define the parameter $A$ as the maximum of any $x$-axis grid boundary $ip$ and its nearest preimage $x^m_j$:$$A = \max_{i} \, \min_{m, j} \left\vert{} x_j^m - ip \right\vert{}.$$ 
Apart from the intervals where $|f'|=1$, we have $\vert{}h_0'\vert{} = 1 + \epsilon_m$ on these central intervals for some small $\epsilon_m$ depending on $\alpha$, hence the preimages $x^m_j$ closely align with the $x$-axis grid boundaries. As $\alpha \to 0$, the derivative of $h_0$ approaches $1$ on all $I_m^0$, causing $A \to 0$. For the intervals where $\vert{}f'\vert{} = 1$, we have $I_m^0 = I_m$ and $\vert{}h_0'\vert{} = 1$. In this situation, the preimages $x_j^m$ coincide with the grid points. This alignment does not modify the forthcoming construction. We select a uniform, global tolerance $q > 0$ to be fixed later and define a buffer distance $B = A + q$. Note that $B \to 0$ as both $\alpha, q \to 0$. For each $x$-axis grid cell $[ip, (i+1)p]$, this global buffer distance defines a central ''bulk'' region $[ip+B, (i+1)p-B]$ bounded by two buffer zones.

We now focus on a specific interval $I_m^0$ corresponding to our chosen $J$ to replace $h_0$ with an oscillating function $h$. For simplicity, assume $h_0$ is increasing on $I_m^0$ (the decreasing case is handled symmetrically).
Because the maximum drift $A \to 0$ as $\alpha \to 0$, each preimage interval $[x_j, x_{j+1}]$ naturally pairs with a unique corresponding grid cell $[ip, (i+1)p]$ if $\alpha$ is sufficiently small.

 We replace \(h_0\) with an oscillating function \(h\) on \(I_m^0\). Between any two consecutive preimages \([x_j, x_{j+1}]\), we define \(h\) piecewise to span the identical vertical range \([jp, (j+1)p]\) as $h_0$ (see Figure \ref{fig6}):
\begin{enumerate}
    \item \textbf{Left transition:} From \(x_j\) (or the left boundary of \(I_m^0\)) to the start of the bulk region \(ip+B\), \(h\) consists of a single ``tent map'' spanning the vertical range \(p\).
    \item \textbf{The bulk:} Between the bounds of the bulk region \([ip+B, (i+1)p-B]\), \(h\) consists of \(L\) identical, consecutive tent maps, each of vertical height \(p\).
    \item \textbf{Right transition:} From the end of the bulk region \((i+1)p-B\) to \(x_{j+1}\), \(h\) forms a single affine slope that reaches the maximum height \((j+1)p\) exactly at \(x_{j+1}\). (This ensures a continuous transition into the next segment, where the pattern repeats).
\end{enumerate}

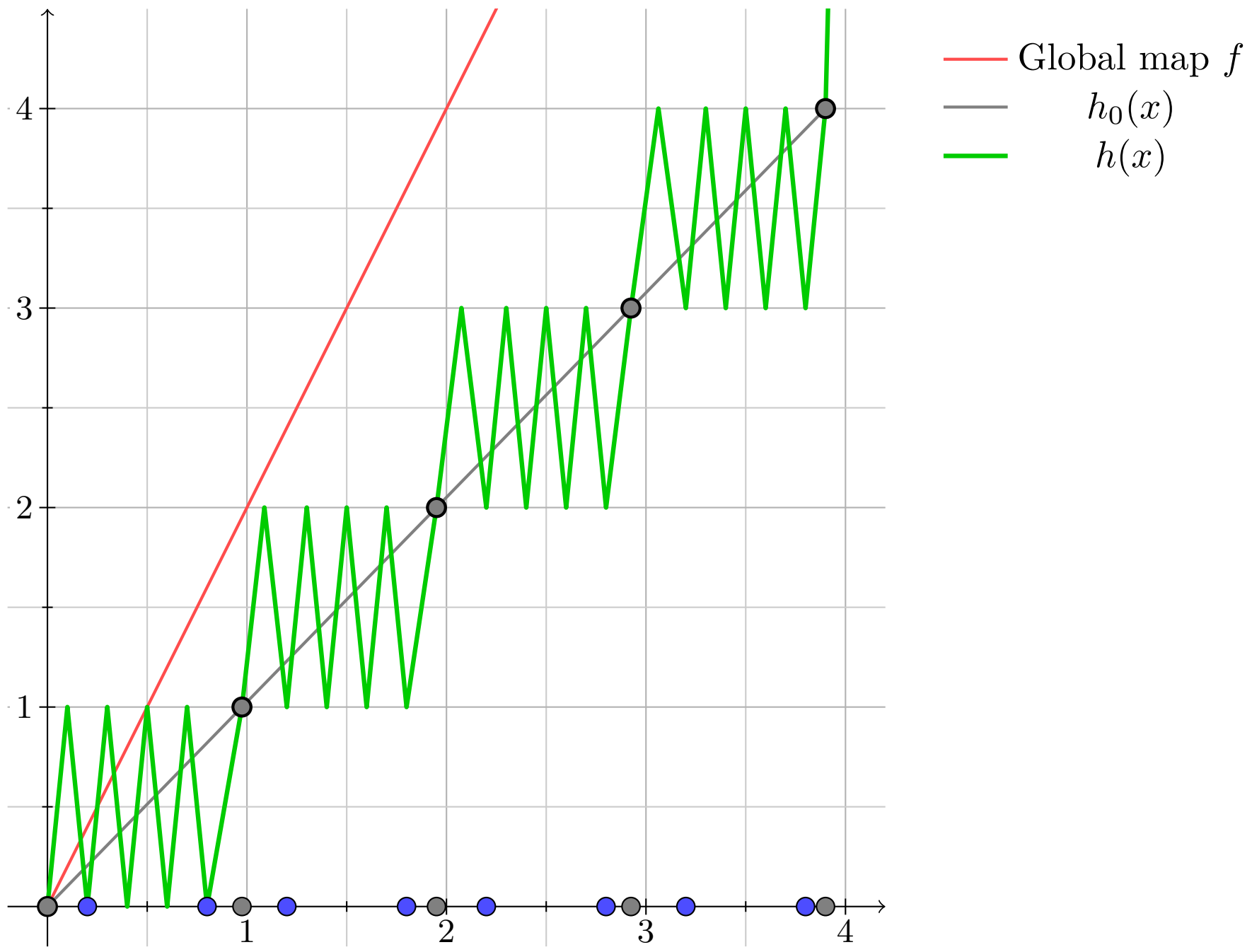
\begin{figure}[htbp]
    \centering
    \begin{tikzpicture}
    \begin{axis}[
        width=12cm,
        unit vector ratio=1 1 1, 
        xmin=-0.2,
        xmax=4.2,
        ymin=-0.2,
        ymax=4.5,
        axis lines=middle,
        axis line style={->},
        xtick={0,1,2,3,4},
        ytick={0,1,2,3,4},
        minor tick num=1, 
        grid=both,
        major grid style={line width=0.4pt, draw=gray!60},
        minor grid style={line width=0.4pt, draw=gray!40},
        tick style={
            draw=black,
            line width=0.4pt,
            preaction={
                draw=white,
                line width=2pt
            }
        },
        tick label style={
            font=\small,
            fill=white,
            inner sep=1.5pt
        },
        legend style={
            at={(1.05, 0.98)},
            anchor=north west,
            draw=none,
            fill=white
        },
        clip=true
    ]

    \addplot[
        thick, 
        red!70
    ] coordinates {
        (0,0) (4,8)
    };
    \addlegendentry{Global map $f$}

    \addplot[thick, gray] coordinates {(-1,-1)};
    \addlegendentry{$h_0(x)$}

    \addplot[
        thick, 
        gray,
        mark=*,
        mark options={fill=gray, draw=black, scale=1.2},
        forget plot 
    ] coordinates {
        (0,0) (0.975, 1) (1.95, 2) (2.925, 3) (3.9, 4)
    };
    
    \addplot[
        thick, 
        gray,
        forget plot
    ] coordinates {
        (3.9, 4) (4, 8)
    };

    \addplot[
        very thick, 
        green!80!black, 
        line join=round
    ] coordinates {
        (0,0) (0.1,1) (0.2,0)
        (0.3,1) (0.4,0) (0.5,1) (0.6,0) (0.7,1) (0.8,0)
        (0.975,1) 
        (1.0875,2) (1.2,1)
        (1.3,2) (1.4,1) (1.5,2) (1.6,1) (1.7,2) (1.8,1)
        (1.95,2)
        (2.075,3) (2.2,2)
        (2.3,3) (2.4,2) (2.5,3) (2.6,2) (2.7,3) (2.8,2)
        (2.925,3)
        (3.0625,4) (3.2,3)
        (3.3,4) (3.4,3) (3.5,4) (3.6,3) (3.7,4) (3.8,3)
        (3.9,4)
        (4,8)
    };
    \addlegendentry{$h(x)$}

    \addplot[
        only marks, 
        mark=*, 
        mark options={fill=gray, draw=black, scale=1.2}
    ] coordinates {
        (0,0) (0.975,0) (1.95,0) (2.925,0) (3.9,0)
    };

    \addplot[
        only marks, 
        mark=*, 
        mark options={fill=blue!70, draw=black, scale=1.2}
    ] coordinates {
        (0.2,0) (0.8,0) 
        (1.2,0) (1.8,0) 
        (2.2,0) (2.8,0) 
        (3.2,0) (3.8,0)
    };

    \end{axis}
    \end{tikzpicture}
    \caption{Construction of the modified function $h$ (green) from the intermediate function $h_0$ (grey). The exact preimages $x_j = j \cdot 0.975$ are marked as grey dots on the $x$-axis. The blue dots indicate the boundaries of the central ``bulk'' regions $[ip+B, (i+1)p-B]$ defined by the buffer distance $B=0.2$. Within each bulk region, $L=3$ identical tent maps are placed, with continuous affine transitions connecting to the preimages $x_j$.}
    \label{fig6}
\end{figure}

In the case where \(h_0\) is decreasing, the local geometry is simply mirrored: we begin with a decreasing slope from \(x_j\) to the bulk, place \(L\) tents in the bulk, and end with a single tent in the right buffer. On the outer intervals \(I_m^-\) and \(I_m^+\), \(h\) retains the identical slope of \(h_0\). 

By construction, the total horizontal domain allocated to traversing any vertical grid segment \([jp, (j+1)p]\) matches exactly the domain \([x_j, x_{j+1}]\) used by \(h_0\). Because \(h\) spans this vertical range using full vertical sweeps, the 1-dimensional Lebesgue measure is perfectly preserved, guaranteeing \(h \in \mathrm{CPA}(\lambda)\). Furthermore, by choosing $L$ to be sufficiently large and $\alpha, q$ to be sufficiently small, we guarantee that $\vert{}h'\vert{} \neq 1$ everywhere, thus \(h \in \mathrm{CPA}^*(\lambda)\).

\vspace{6pt}
We show that by choosing the parameters $\alpha, q, L$ of the structure above appropriately, we can achieve $\Hp_k(h)>(1-\delta)k\log k$ for some $k>N$.

First, let us observe the macroscopic structure of $h$. Recall that the initial approximation $g$ is a piecewise affine, measure-preserving map with slopes $\pm 1$, and its partition intervals align perfectly with our uniform $p \times p$ grid. Due to the measure-preserving property of $f$, the graph of $g$ occupies the grid exactly like a permutation matrix: every horizontal row and vertical column of the grid contains exactly one diagonal segment of $g$. While our subsequent modifications slightly shifted the horizontal boundaries from the rigid grid points $ip$ to the preimages $x_j$, the final map $h$ still maps the domain $[x_j, x_{j+1}]$ strictly onto the vertical range $[jp, (j+1)p]$. Thus, the global dynamics of $h$ preserve this permutation matrix-like routing on the ''bulk'' regions. 

Inside each of these intervals, $h$ shares an identical micro-structure. In the bulk region, it consists of $L$ consecutive tent maps, which provide $2L$ strictly monotone branches, each mapping onto the full vertical range $p$. We can analyze the resulting complexity using a standard symbolic dynamics approach. We assign labels $1, 2, \dots, 2L$ to these small branches comprising the bulk region. If the first $k$ iterates of a point visit branches with $k$ distinct labels, its symbolic trajectory uniquely dictates the relative vertical ordering of its orbit, and thus completely determines its length-$k$ permutation pattern.

Because each of these $2L$ branches spans the full height of the cell, every $k$-length permutation is realized on a subset whose measure is roughly proportional to the Lebesgue measure of the bulk region divided by $k!$. By leveraging the periodic macro-structure of $h$, we will show that we observe enough distinct permutation patterns with sufficiently uniform probabilities such that the permutation entropy $\Hp_k(h)$ is on the order of $k \log k$. The rigorous details of these estimates are carried out below.

\vspace{6pt}
First, we establish a lower bound for the measure of each $k$-length permutation by analyzing the micro-dynamics within a single grid cell. To simplify the calculation, we view this local behavior as being scaled to a $[0,1] \to [0,1]$ map, which leaves a bulk region of relative measure $1-2B/p$. 

Since each of the $2L$ strictly monotone branches in the bulk region spans the entire vertical range in a linear manner, the transition from any branch $i$ to any specific branch $j$ in the subsequent cell occurs with probability $\frac{1}{2L}\left(1-\frac{2B}{p}\right)$. 

This uniform transition property implies that the relative measure of points whose first $k$ iterates remain entirely within the bulk regions and form a symbolic sequence of $k$ distinct branches is exactly:
\[ \gamma :=\left(1 - \frac{2B}{p}\right)^k \frac{2L \cdot (2L-1) \cdots (2L-k+1)}{(2L)^k}. \]
By choosing our parameters such that $B \to 0$ and $L \to \infty$, this proportion tends to 1. This guarantees that we observe these ``good'' symbolic sequence—trajectories composed of $k$ different branch entries—on a subset whose measure is arbitrarily close to the full measure of the cell. 

Crucially, because the vertical ordering of the orbit is uniquely determined by the chosen branches, every $k$-length permutation is generated by exactly $\binom{2L}{k}$ of these good symbolic series. Since all such series share the exact same measure due to the uniformity of the branches, it follows that we obtain all $k$-length permutations on this large bulk set with uniform measures.

We can now estimate the permutation entropy of our global function $h$ using Lemma \ref{2}.
The permutation matrix structure decomposes the grid into disjoint cycles. Let us take a cycle of length $j$ ($1 \le j \le 1/p$). On a macroscale, the sequence of $k$ iterates of any point in this cycle is distributed periodically among the $j$ grid cells. 

Because the grid separates the cells into disjoint vertical intervals, the relative ordering between any two iterates that land in different cells is strictly fixed by the macroscopic cycle. Thus, any global length-$k$ permutation is uniquely determined simply by the internal relative ordering of the iterates that land within the same cell. 

Recall that for the original $L$-tent map, all $k!$ permutations appear with equal measure on our designated ``good'' set. To rigorously capture why the local permutations multiply, we can view a trajectory through a two-coordinate lens: the first coordinate tracks the macroscopic movement (the cycle through the $j$ grid cells), while the second coordinate tracks the microscopic movement (the choice among the $2L$ affine branches within a cell). Because the map expands linearly over the full vertical extent of the chosen cells, the choices of the second coordinate (the specific branches taken) are entirely independent as the orbit cycles through the $j$ cells. This ensures that the relative orderings of the iterates within distinct cells can be formed completely freely without restricting one another. Over $k$ iterates, the periodic cycle ensures that the orbit visits the $j$ cells along the cycle $t_1, \dots, t_j$ times, respectively, where each $t_i\in \{\lfloor k/j \rfloor, \lceil k/j\rceil\}$. The visits to any specific cell correspond to time steps belonging to a single modulo $j$ residue class (i.e., taking every $j$-th step). We can determine the distribution of these local patterns directly from the inherent combinatorial symmetries of the permutations. If we take the uniformly distributed $k!$ permutations of the original map and extract the coordinates at every $j$-th step, symmetry guarantees that in the $i$th cell, every possible relative ordering of length $t_i$ appears exactly the same number of times. Thus, we observe all $t_i !$ orderings with the same frequency within each individual cell. Moreover, as the underlying length $k$ permutation was uniform random, these $j$ many permutations of length $t_i$ are independent from each other.
This yields all $\prod_{i=1}^{j}t_i$ combined permutations on subsets of identical measure.  Furthermore, starting from different grid cells in the cycle yields distinct macro-routes, providing an additional factor of $j$, as $k>1/p\geq j$. Thus, across the union of these good sets, we observe at least $j \cdot \prod_{i=1}^{j}t_i$ distinct permutations, all distributed uniformly in measure. The entropy of this can be estimated from below by replacing each $t_i$ by $\lfloor k/j \rfloor$.

Using Lemma \ref{2}, the permutation entropy of the entire map $h$ is bounded below by the convex combination of the entropies of its disjoint cycles. Letting $\alpha_i$ denote the weights of the initial macro-cycle blocks and $j_i$ its length, we obtain the global lower bound on the normalized entropy: 
\[ \frac{\Hp_k(h)}{k \log k} \ge \sum_i \alpha_i \left[ \gamma \frac{\log \big( j_i \cdot (\lfloor k/j_i \rfloor !)^{j_i} \big)}{k \log k} \right]. \] 

Taking $k\to \infty$, by Stirling's approximation, the innermost fractions go to 1 for each $i$. Thus we can fix $k$ such that each of these fractions is at least $(1-\delta)^{1/2}$. Given this choice, we can fix small $\alpha, q$ (and through that, small enough $B$) and large enough $L$ such that $\gamma\geq(1-\delta)^{1/2}$.
Since $\sum_i \alpha_i = 1$, with such choices we recover
$$\frac{\Hp_k(h)}{k \log k}\geq 1-\delta.$$
This concludes the proof.
\end{proof}

\section*{Declaration of AI usage}

During the preparation of this document, Google Gemini was used to improve the writing of the original proofs, assist with LaTeX formatting, and generate TikZ code for figures. ChatGPT (GPT-6) was used for proofreading the essentially finished manuscript. The authors take full responsibility for the content of this work.

\printbibliography

\end{document}